# A Multivariate Weibull Disitribution

Cheng K. Lee
chengli@uab.edu
*Charlotte, North Carolina, USA*

Miin-Jye Wen
*National Cheng Kung University, City Tainan, Taiwan, R.O.C.*

**Summary.**

A multivariate survival function of Weibull Distribution is developed by expanding the theorem by Lu and Bhattacharyya (1990). From the survival function, the probability density function, the cumulative probability function, the determinant of the Jacobian Matrix, and the general moment are derived. The proposed model is also applied to the tumor appearance data of female rats.

## 1. Introduction

Lu and Bhattacharyya (1990) developed a joint survival function by letting $h_1(x)$ and $h_2(y)$ be two arbitrary failure rate functions on $[0,\infty)$, and $H_1(x)$ and $H_2(y)$ their corresponding cumulative failure rate. Given the stress $S=s > 0$, the joint survival function conditioned on $s$, as they defined, is

$$\overline{F}(x,y\,|\,s) = \exp\left\{-\left[H_1(x)+H_2(y)\right]^{\gamma} s\right\},$$

where $\gamma$ measures the conditional association of $X$ and $Y$. Further, based on the joint survival function, they proved a theorem that a bivariate survival function $\overline{F}(x,y\,|\,s)$ can be derived with the marginals $\overline{F}_x$ and $\overline{F}_y$ given the assumption that the Laplace transform of the stress $S$ exists on $[0,\infty)$ and is strictly decreasing.
From the theorem, they derived a bivariate Weibull Distribution

$$\overline{F}(x,y) = \exp\left\{-\left[\left(\frac{x}{\lambda_1}\right)^{\frac{\gamma_1}{\alpha}} + \left(\frac{y}{\lambda_2}\right)^{\frac{\gamma_2}{\alpha}}\right]^{\alpha}\right\},$$

where $0 < \alpha \leq 1$, $0 < \lambda_1, \lambda_2 < \infty$, and $0 < \gamma_1, \gamma_2 < \infty$. This bivariate Weibull Distribution is exactly the same as developed by Hougaard (1986).

By the same steps, the theorem can be expanded to more than two random variables, and, therefore, a multivariate survival function of Weibull distribution is constructed as

$$S(x_1, x_2, \ldots, x_n) = \exp\left\{-\left[\left(\frac{x_1}{\lambda_1}\right)^{\frac{\gamma_1}{\alpha}} + \left(\frac{x_2}{\lambda_2}\right)^{\frac{\gamma_2}{\alpha}} + \ldots + \left(\frac{x_n}{\lambda_n}\right)^{\frac{\gamma_n}{\alpha}}\right]^{\alpha}\right\}, \tag{1}$$

where $\alpha$ measures the association among the variables, $0 < \alpha \leq 1$, $0 < \lambda_1, \lambda_2, \ldots, \lambda_n < \infty$, and $0 < \gamma_1, \gamma_2, \ldots, \gamma_n < \infty$.

This model can also be derived by using a copula construction in which the generator is $(-\log(.))^{\alpha}$ (Frees and Valdez, 1998). Equation 1 is similar to the "genuine multivariate Weibull distribution" developed by Crowder (1989) who studied another version extended from the genuine multivariate Weibull distribution.

In this paper, we mathematically intensively studied the proposed multivariate Weibull model of Equation 1 by giving the probability density function in section 2, the Jacobian matrix in section 3, the general moment in section 4 and an application in section 5.

## 2. Probability Density Function of The Multivariate Weibull Distribution

The multivariate probability density function $f(x_1, x_2, \ldots, x_n)$ of a multivariate distribution function can be obtained by differentiating the multivariate survival function with respect to each variable. Li (1997) has shown that

$$f(x_1, x_2, \ldots, x_n) = (-1)^n \frac{\partial^n S(x_1, x_2, \ldots, x_n)}{\partial x_1 \partial x_2 \ldots \partial x_n}.$$

Using Li's derivation and one of the special cases of the multivariate Faa di Bruno formula by Constantine and Savits (1996), the probability density function is

$$f(x_1, x_2, \ldots, x_n) = \left(\frac{-1}{\alpha}\right)^n \exp\left\{-\left[\left(\frac{x_1}{\lambda_1}\right)^{\frac{\gamma_1}{\alpha}} + \left(\frac{x_2}{\lambda_2}\right)^{\frac{\gamma_2}{\alpha}} + \ldots + \left(\frac{x_n}{\lambda_n}\right)^{\frac{\gamma_n}{\alpha}}\right]^{\alpha}\right\}$$

$$\cdot \left[\left(\frac{\gamma_1}{\lambda_1}\right)\left(\frac{\gamma_2}{\lambda_2}\right)\cdots\left(\frac{\gamma_n}{\lambda_n}\right)\right]\left[\left(\frac{x_1}{\lambda_1}\right)^{\frac{\gamma_1}{\alpha}-1}\left(\frac{x_2}{\lambda_2}\right)^{\frac{\gamma_2}{\alpha}-1}\cdots\left(\frac{x_n}{\lambda_n}\right)^{\frac{\gamma_n}{\alpha}-1}\right]$$

$$\cdot \sum_{i=1}^{P(n)} \left\{(-1)^{k_i} P_s(n,i)\left(\prod_{j=1}^{k_i} \alpha^{\underline{n_j}}\right)\left[\left(\frac{x_1}{\lambda_1}\right)^{\frac{\gamma_1}{\alpha}} + \left(\frac{x_2}{\lambda_2}\right)^{\frac{\gamma_2}{\alpha}} + \ldots + \left(\frac{x_n}{\lambda_n}\right)^{\frac{\gamma_n}{\alpha}}\right]^{k_i \alpha - n}\right\}, \tag{2}$$

where $k_i$ is the number of summands of the $i$th partition of $n$ such that $n_1 + n_2 + \cdots + n_{k_i} = n$, $n_1 \geq n_2 \geq \cdots \geq n_{k_i} > 0$, $1 \leq k_i \leq n$; $\alpha^{\underline{n_j}}$ is equal to $\alpha(\alpha-1)\ldots(\alpha-n_j+1)$, the falling factorial of $\alpha$ (Kunth, 1992); $P(n)$ is the total number of partitions of $n$; $P_s(n,i)$ is the total number of set partitions of the set $S_n = \{1,\ldots,n\}$ corresponding to the $i$th partition of $n$. The specific way of partitioning $n$ and $S_n$ is given by McCullagh and Wilks (1988). In their paper, partitions of $n$ are in increasing number of summands and ordering all the summands in inverse lexicographic order when a partition has the same number of summands, and $S_n = \{1, \ldots, n\}$ is partitioned by "listing the blocks from the largest to the smallest and by breaking the ties of equal sized blocks by ordering them lexicographically" and the number of blocks in a set partition is equal to the number of summands of the corresponding partition of $n$. For example, the total number of blocks of the partition of $S_n$ corresponding to the $i$th partition is $k_i$ and the numbers of elements in each block are equal to $n_1, n_2, \cdots, n_{k_i}$.

## 3. The Jacobian Matrix

Similar to the derivation of the Bivariate Weibull Distribution by Lu and Bhattacharyya (1990), let $(y_1, y_2, \ldots, y_n)$

$$= \left( \frac{z_1}{z_1 + z_2 + \ldots + z_n}, \frac{z_2}{z_1 + z_2 + \ldots + z_n}, \ldots, \frac{z_{n-1}}{z_1 + z_2 + \ldots + z_n}, (z_1 + z_2 + \ldots + z_n)^\alpha \right) \quad (3)$$

where $z_1 = \left(\frac{x_1}{\lambda_1}\right)^{\frac{\gamma_1}{\alpha}}, z_2 = \left(\frac{x_2}{\lambda_2}\right)^{\frac{\gamma_2}{\alpha}}, \ldots, z_n = \left(\frac{x_n}{\lambda_n}\right)^{\frac{\gamma_n}{\alpha}}$.

Then, $x_1 = y_1^{\frac{\alpha}{\gamma_1}} y_n^{\frac{1}{\gamma_1}} \lambda_1, x_2 = y_2^{\frac{\alpha}{\gamma_2}} y_n^{\frac{1}{\gamma_2}} \lambda_2, \ldots, x_{n-1} = y_{n-1}^{\frac{\alpha}{\gamma_{n-1}}} y_n^{\frac{1}{\gamma_{n-1}}} \lambda_{n-1}$,

$x_n = (1 - y_1 - y_2 - \cdots - y_{n-1})^{\frac{\alpha}{\gamma_n}} y_n^{\frac{1}{\gamma_n}} \lambda_n$.

Note that $z_1, z_2, \ldots, z_n > 0$, and
$1 - y_1 - y_2 - \ldots - y_{n-1}$
$= 1 - \frac{z_1 + z_2 + \ldots + z_{n-1}}{z_1 + z_2 + \ldots + z_n}$
$= \frac{z_n}{z_1 + z_2 + \ldots + z_n} > 0$.

The Jacobian matrix is

$$J = \begin{vmatrix} \frac{\partial x_1}{\partial y_1} & \frac{\partial x_1}{\partial y_2} & \cdots & \frac{\partial x_1}{\partial y_n} \\ \frac{\partial x_2}{\partial y_1} & \frac{\partial x_2}{\partial y_2} & \cdots & \frac{\partial x_2}{\partial y_n} \\ \vdots & \vdots & & \vdots \\ \frac{\partial x_n}{\partial y_1} & \frac{\partial x_n}{\partial y_2} & \cdots & \frac{\partial x_n}{\partial y_n} \end{vmatrix}.$$

Let $C(i,j)$ be the $i$th row and $j$th column in the Jacobian matrix, then

$$C(i,i) = \frac{\alpha \lambda_i y_i^{\frac{\alpha}{\gamma_i}-1} y_n^{\frac{1}{\gamma_i}}}{\gamma_i}, \quad i=1,2,\ldots,n-1,$$

$$C(i,n) = \frac{\lambda_i y_i^{\frac{\alpha}{\gamma_i}} y_n^{\frac{1}{\gamma_i}-1}}{\gamma_i}, \quad i=1,2,\ldots,n-1,$$

$$C(n,i) = -\frac{\alpha \lambda_n (1-y_1-y_2-\cdots-y_{n-1})^{\frac{\alpha}{\gamma_n}-1} y_n^{\frac{1}{\gamma_n}}}{\gamma_n}, \quad i=1,2,\ldots,n-1,$$

$$C(n,n) = \frac{\lambda_n (1-y_1-y_2-\cdots-y_{n-1})^{\frac{\alpha}{\gamma_n}} y_n^{\frac{1}{\gamma_n}-1}}{\gamma_n},$$

$C(i,j)=0$ when $i \neq j$, $i=2,\ldots,n-1$, $j=2,\ldots,n-1$.

The determinant of the Jacobian matrix can be obtained using Gaussian elimination to construct an upper triangle matrix.
Then, the determinant is equal to the product of the diagonal elements.

$$|J| = \left(\prod_{i=1}^{n-1} C(i,i)\right)\left(C(n,n) - \sum_{j=1}^{n-1} \frac{C(n,j)C(j,n)}{C(j,j)}\right)$$

$$= \frac{\alpha^{n-1} \lambda_1 \lambda_2 \cdots \lambda_n y_1^{\frac{\alpha}{\gamma_1}-1} y_2^{\frac{\alpha}{\gamma_2}-1} \cdots y_{n-1}^{\frac{\alpha}{\gamma_{n-1}}-1} (1-y_1-y_2-\cdots-y_{n-1})^{\frac{\alpha}{\gamma_n}-1} y_n^{\frac{1}{\gamma_1}+\frac{1}{\gamma_2}+\cdots+\frac{1}{\gamma_n}-1}}{\gamma_1 \gamma_2 \cdots \gamma_n}. \quad (4)$$

After the derivation of the Jacobian, the PDF in terms of $y_1, y_2, \cdots, y_n$,

$$g(y_1, y_2, \cdots, y_n) =$$

$$f\left(y_1^{\frac{\alpha}{\gamma_1}} y_n^{\frac{1}{\gamma_1}} \lambda_1, y_2^{\frac{\alpha}{\gamma_2}} y_n^{\frac{1}{\gamma_2}} \lambda_2, \cdots, y_{n-1}^{\frac{\alpha}{\gamma_{n-1}}} y_n^{\frac{1}{\gamma_{n-1}}} \lambda_{n-1}, (1-y_1-y_2-\cdots-y_{n-1})^{\frac{\alpha}{\gamma_n}} y_n^{\frac{1}{\gamma_n}} \lambda_n\right)|J|$$

$$= (-1)^n \left\{\sum_{i=1}^{P(n)} \left(\alpha^{-1} y_n^{k_i-1} (-1)^{k_i} P_s(n,i) \left(\prod_{j=1}^{k_i} \alpha^{\frac{n_j}{}}\right)\right)\right\} \exp(-y_n)$$

$$= \left(\Gamma(n) y_1^{1-1} y_2^{1-1} \cdots y_{n-1}^{1-1} (1-y_1-y_2-\cdots-y_{n-1})^{1-1}\right) f(y_n), \quad (5)$$

where $y_1, y_2, \ldots, y_{n-1}$ has a Dirichlet distribution with the probability density equal to $\Gamma(n) y_1^{1-1} y_2^{1-1} \cdots y_{n-1}^{1-1} (1-y_1-y_2-\cdots-y_{n-1})^{1-1}$, and

$$f(y_n) = \frac{(-1)^n}{\Gamma(n)} \left\{\sum_{i=1}^{P(n)} \left(\alpha^{-1} y_n^{k_i-1} (-1)^{k_i} P_s(n,i) \left(\prod_{j=1}^{k_i} \alpha^{\frac{n_j}{}}\right)\right)\right\} \exp(-y_n) \quad (6)$$

has a mixture distribution of the exponential distribution and Gamma distribution.

Equation (6) can be rewritten as

$$f(y_n) = \frac{(-1)^n}{\Gamma(n)} \left\{ \sum_{i=1}^{P(n)} \left( \alpha^{-1} y_n^{k_i-1} (-1)^{k_i} \Gamma(k_i) P_s(n,i) \left( \prod_{j=1}^{k_i} \alpha^{n_j} \right) \right) \frac{\exp(-y_n)}{\Gamma(k_i)} \right\}.$$

When it is integrated over the range of $y_n$, it becomes

$$\frac{(-1)^n}{\Gamma(n)} \left\{ \sum_{i=1}^{P(n)} \left( \alpha^{-1} (-1)^{k_i} \Gamma(k_i) P_s(n,i) \left( \prod_{j=1}^{k_i} \alpha^{n_j} \right) \right) \right\} = 1.$$

That is the weights of $y_n$ are summed to 1. The probability density function of $y_n$ is the mixed Gamma distribution by Downton (1969). Following his derivation, the cumulative density function of $y_n$ is

$$1 - \frac{(-1)^n}{\Gamma(n)} \left( \sum_{i=1}^{n} \frac{y_n^{i-1}}{\Gamma(i)} \sum_{k \geq i}^{P(n)} \left( \alpha^{-1} (-1)^k \Gamma(k) P_s(n,i) \left( \prod_{j=1}^{k} \alpha^{n_j} \right) \right) \right) \exp(-y_n). \qquad (7)$$

### 4. The General Moment

The general moment of $x_1, x_2, \cdots, x_n$ is

$$E\left[ x_1^{i_1} x_2^{i_2} \cdots x_n^{i_n} \right] =$$

$$E\left[ \left( y_1^{\frac{\alpha}{\gamma_1}} y_n^{\frac{1}{\gamma_1}} \lambda_1 \right)^{i_1} \left( y_2^{\frac{\alpha}{\gamma_2}} y_n^{\frac{1}{\gamma_2}} \lambda_2 \right)^{i_2} \cdots \left( y_{n-1}^{\frac{\alpha}{\gamma_{n-1}}} y_n^{\frac{1}{\gamma_{n-1}}} \lambda_{n-1} \right)^{i_{n-1}} \left( (1 - y_1 - y_2 - \cdots - y_{n-1})^{\frac{\alpha}{\gamma_n}} y_n^{\frac{1}{\gamma_n}} \lambda_n \right)^{i_n} \right]$$

$$= \lambda_1^{i_1} \lambda_2^{i_2} \cdots \lambda_n^{i_n} E\left[ y_1^{\frac{i_1 \alpha}{\gamma_1}} y_2^{\frac{i_2 \alpha}{\gamma_2}} \cdots y_{n-1}^{\frac{i_{n-1} \alpha}{\gamma_{n-1}}} (1 - y_1 - y_2 - \cdots - y_{n-1})^{\frac{i_n \alpha}{\gamma_n}} \right] E\left[ y_n^{\frac{i_1}{\gamma_1} + \frac{i_2}{\gamma_2} + \cdots \frac{i_n}{\gamma_n}} \right]$$

Then,

$$E\left[ y_1^{\frac{i_1 \alpha}{\gamma_1}} y_2^{\frac{i_2 \alpha}{\gamma_2}} \cdots y_{n-1}^{\frac{i_{n-1} \alpha}{\gamma_{n-1}}} (1 - y_1 - y_2 - \cdots - y_{n-1})^{\frac{i_n \alpha}{\gamma_n}} \right]$$

$$= \Gamma(n) \int \cdots \int y_1^{\frac{i_1 \alpha}{\gamma_1}} y_2^{\frac{i_2 \alpha}{\gamma_2}} \cdots y_{n-1}^{\frac{i_{n-1} \alpha}{\gamma_{n-1}}} (1 - y_1 - y_2 - \cdots - y_{n-1})^{\frac{i_n \alpha}{\gamma_n}} dy_1 \cdots dy_{n-1}$$

$$= \frac{\Gamma(n) \Gamma\left( \frac{i_1 \alpha}{r_1} + 1 \right) \Gamma\left( \frac{i_2 \alpha}{r_2} + 1 \right) \cdots \Gamma\left( \frac{i_n \alpha}{r_n} + 1 \right)}{\Gamma\left( \alpha \left( \frac{i_1}{\gamma_1} + \frac{i_2}{\gamma_2} + \cdots + \frac{i_n}{\gamma_n} \right) + n \right)}$$

which is the Dirichlet integral (Rao, 1954.)

For $E\left[ y_n^{\frac{i_1}{\gamma_1} + \frac{i_2}{\gamma_2} + \cdots \frac{i_n}{\gamma_n}} \right]$, let $\frac{i_1}{\gamma_1} + \frac{i_2}{\gamma_2} + \cdots + \frac{i_n}{\gamma_n} = c$, then

$$E\left[ y_n^c \right]$$

$$= \frac{(-1)^n}{\Gamma(n)} \int_0^\infty y_n^c \sum_{i=1}^{P(n)} \left\{ \alpha^{-1} y_n^{k_i-1} (-1)^{k_i} P_s(n,i) \left( \prod_{j=1}^{k_i} \alpha^{\underline{n_j}} \right) \right\} \exp(-y_n) dy_n$$

$$= \frac{(-1)^n}{\Gamma(n)} \sum_{i=1}^{P(n)} \left\{ \left[ \alpha^{-1} (-1)^{k_i} P_s(n,i) \left( \prod_{j=1}^{k_i} \alpha^{\underline{n_j}} \right) \right] \int_0^\infty y_n^{c+k_i-1} \exp(-y_n) dy_n \right\}$$

$$= \frac{(-1)^n}{\Gamma(n)} \sum_{i=1}^{P(n)} \left[ \alpha^{-1} (-1)^{k_i} P_s(n,i) \left( \prod_{j=1}^{k_i} \alpha^{\underline{n_j}} \right) \Gamma(c+k_i) \right]$$

$$= \frac{(-1)^n}{\Gamma(n)} \sum_{i=1}^{P(n)} \left[ \alpha^{-1} (-1)^{k_i} P_s(n,i) \left( \prod_{j=1}^{k_i} \alpha^{\underline{n_j}} \right) \Gamma(c) c^{\overline{k_i}} \right] \tag{8}$$

where $c^{\overline{k_i}}$ is the rising factorial defined as $c(c+1)\cdots(c+k_i-1)$ by Knuth (1992).

Considering $P_s(n,i)$ in the above equation, it is the total number of set partitions corresponding to the $i$th partition of $n$ such that $n_1 + n_2 + \cdots + n_{k_i} = n$, $n_1, n_2, \ldots, n_{k_i} > 0$. It has been shown by McCullaph and Wilks (1988) that

$$P_s(n,i) = \frac{n!}{n_1! n_2! \cdots n_{k_i}! m_1! m_2! \cdots m_d!}$$

where $m_1, m_2, \ldots, m_d$ are the number of each distinct summand.

Then, the product $P_s(n,i) \left( \prod_{j=1}^{k_i} \alpha^{\underline{n_j}} \right)$

$$= \frac{n!}{n_1! n_2! \cdots n_{k_i}! m_1! m_2! \cdots m_d!} \alpha^{\underline{n_1}} \alpha^{\underline{n_2}} \cdots \alpha^{\underline{n_{k_i}}}$$

$$= \frac{n!}{m_1! m_2! \cdots m_d!} \frac{\alpha!}{n_1!(\alpha-n_1)!} \frac{\alpha!}{n_2!(\alpha-n_2)!} \cdots \frac{\alpha!}{n_{k_i}!(\alpha-n_{k_i})!}$$

$$= \frac{n!}{m_1! m_2! \cdots m_d!} \binom{\alpha}{n_1} \binom{\alpha}{n_2} \cdots \binom{\alpha}{n_{k_i}}$$

$$= \frac{n!}{k_i!} \frac{k_i!}{m_1! m_2! \cdots m_d!} \binom{\alpha}{n_1} \binom{\alpha}{n_2} \cdots \binom{\alpha}{n_{k_i}},$$

where $\dfrac{k_i!}{m_1! m_2! \cdots m_d!}$ is the number of permutations of $n_1, n_2, \cdots, n_{k_i}$ of every possible order.

When sum $P_s(n,i) \left( \prod_{j=1}^{k_i} \alpha^{\underline{n_j}} \right)$ over the same value of $k_i$, say, $k$, then,

$$\sum_{k_i=k} P_s(n,i) \left( \prod_{j=1}^{k_i} \alpha^{\underline{n_j}} \right)$$

$$= \sum_{k_i = k} \left[ \frac{n!}{k_i! \, m_1! m_2! \cdots m_d!} \frac{k_i!}{} \binom{\alpha}{n_1} \binom{\alpha}{n_2} \cdots \binom{\alpha}{n_{k_i}} \right]$$

$$= \frac{n!}{k!} \sum_{k_i = k} \left[ \frac{k!}{m_1! m_2! \cdots m_d!} \binom{\alpha}{n_1} \binom{\alpha}{n_2} \cdots \binom{\alpha}{n_k} \right]$$

which equals $C(n, k, \alpha)$, the C-numbers defined by Charalambides (1977). Note that the summation is over all the permutations of $n_{k_i}$ with $k_i = k$.

Using the equality $(-1)^{k_i} c^{\overline{k_i}} = (-c)^{\underline{k_i}}$ (Goldman, Joichi, Reiner and White, 1976),

$$E\left[y_n^c\right]$$

$$= \frac{(-1)^n}{\Gamma(n)} \sum_{i=1}^{P(n)} \left[ \alpha^{-1} (-1)^{k_i} P_s(n, i) \left( \prod_{j=1}^{k_i} \alpha^{\underline{n_j}} \right) \Gamma(c) c^{\overline{k_i}} \right]$$

$$= \frac{(-1)^n}{\Gamma(n)} \alpha^{-1} \Gamma(c) \left\{ \sum_{i=1}^{P(n)} \left[ P_s(n, i) \left( \prod_{j=1}^{k_i} \alpha^{\underline{n_j}} \right) (-c)^{\underline{k_i}} \right] \right\}$$

$$= \frac{(-1)^n}{\Gamma(n)} \alpha^{-1} \Gamma(c) \left\{ \sum_{k=1}^{n} \left[ \sum_{k_i = k} \left( P_s(n, i) \left( \prod_{j=1}^{k_i} \alpha^{\underline{n_j}} \right) (-c)^{\underline{k_i}} \right) \right] \right\}$$

$$= \frac{(-1)^n}{\Gamma(n)} \alpha^{-1} \Gamma(c) \sum_{k=1}^{n} \left[ C(n, k, \alpha) (-c)^{\underline{k}} \right]$$

$$= \frac{(-1)^n}{\Gamma(n)} \alpha^{-1} \Gamma(c) (-\alpha c)^{\underline{n}} \quad \text{(using equation 1.3 by Charalambides, 1977)}$$

$$= \frac{(-1)^n}{\Gamma(n)} \alpha^{-1} \Gamma(c) (-1)^n (\alpha c)^{\overline{n}} \quad \text{(using the formula by Goldman } et\ al.\ 1976\text{)}$$

$$= \frac{1}{\Gamma(n)} (\alpha c + 1)(\alpha c + 2) \cdots (\alpha c + n - 1) \Gamma(c + 1).$$

Therefore, the general moment of $x_1, x_2, \cdots, x_n$ is

$$E\left[ x_1^{i_1} x_2^{i_2} \cdots x_n^{i_n} \right]$$

$$= \lambda_1^{i_1} \lambda_2^{i_2} \cdots \lambda_n^{i_n} E\left[ y_1^{\frac{i_1 \alpha}{\gamma_1}} y_2^{\frac{i_2 \alpha}{\gamma_2}} \cdots y_{n-1}^{\frac{i_{n-1} \alpha}{\gamma_{n-1}}} (1 - y_1 - y_2 - \cdots - y_{n-1})^{\frac{i_n \alpha}{\gamma_n}} \right] E\left[ y_n^{\frac{i_1}{\gamma_1} + \frac{i_2}{\gamma_2} + \cdots \frac{i_n}{\gamma_n}} \right]$$

$$= \lambda_1^{i_1} \lambda_2^{i_2} \cdots \lambda_n^{i_n} \frac{\Gamma(n) \Gamma\left( \frac{i_1 \alpha}{r_1} + 1 \right) \Gamma\left( \frac{i_2 \alpha}{r_2} + 1 \right) \cdots \Gamma\left( \frac{i_n \alpha}{r_n} + 1 \right)}{\Gamma\left[ \alpha \left( \frac{i_1}{\gamma_1} + \frac{i_2}{\gamma_2} + \cdots + \frac{i_n}{\gamma_n} \right) + n \right]}$$

$$\cdot \frac{1}{\Gamma(n)} \left[ \alpha \left( \frac{i_1}{\gamma_1} + \frac{i_2}{\gamma_2} + \cdots + \frac{i_n}{\gamma_n} \right) + 1 \right] \left[ \alpha \left( \frac{i_1}{\gamma_1} + \frac{i_2}{\gamma_2} + \cdots + \frac{i_n}{\gamma_n} \right) + 2 \right] \cdots$$

$$\left[ \alpha \left( \frac{i_1}{\gamma_1} + \frac{i_2}{\gamma_2} + \cdots + \frac{i_n}{\gamma_n} \right) + (n-1) \right] \Gamma \left[ \alpha \left( \frac{i_1}{\gamma_1} + \frac{i_2}{\gamma_2} + \cdots + \frac{i_n}{\gamma_n} \right) + 1 \right]$$

$$= \frac{\lambda_1^{i_1} \lambda_2^{i_2} \cdots \lambda_n^{i_n} \Gamma\left(\frac{i_1 \alpha}{r_1} + 1\right) \Gamma\left(\frac{i_2 \alpha}{r_2} + 1\right) \cdots \Gamma\left(\frac{i_n \alpha}{r_n} + 1\right) \Gamma\left[\left(\frac{i_1}{\gamma_1} + \frac{i_2}{\gamma_2} + \cdots + \frac{i_n}{\gamma_n}\right) + 1\right]}{\Gamma\left[\alpha\left(\frac{i_1}{\gamma_1} + \frac{i_2}{\gamma_2} + \cdots + \frac{i_n}{\gamma_n}\right) + 1\right]} \quad (9)$$

From the general moment, the expectation and the variance of any random variable, and the covariance and the correlation coefficient of any number of random variables can be derived.

**5. Application**
We analyze the data published by Mantel, Bohidar and Ciminear (1997). The data (Table 1) contains 50 litters of female rats with one drug-treated and two control rats in each letter. The same data was also analyses by Hougaard (1986) using a bivariate Weibull distribution. We assume the time to the appearance of tumor of the treatment group, the control group 1 and the control group 2 are Weibull distributed. Table 2 displays the counts of combinations of censoring status of litters. The results of parameter estimates and standard errors based on the second derivatives valued at the maximized log-likelihood function are in Table 3. The estimates for the 3 shape parameters are significantly greater than 1 with significant level equal to 0.05 indicating that the 3 group have a monotonically increasing hazard function. The estimate of the association parameter $\alpha$ is not significantly different from 1 indicating that the time to the tumor occurrence among the 3 groups are not associated. Without considering the standard errors, by equation (9), the correlation coefficient of the treatment group and control group 1 is 0.163. The correlation coefficient of the treatment group and the control group 2 is 0.159. The correlation coefficient of the two control groups is 0.157.

Table 1. Time to tumor appearance in weeks of treatment group (T) and control groups (C1, C2).

| Litter | T | C1 | C2 | Litter | T | C1 | C2 |
|---|---|---|---|---|---|---|---|
| 1 | 101+ | 49+ | 104+ | 26 | 104+ | 102+ | 104+ |
| 2 | 104+ | 104+ | 104+ | 27 | 77+ | 97+ | 79+ |
| 3 | 89+ | 104+ | 104+ | 28 | 88 | 96 | 104 |
| 4 | 104 | 94 | 77 | 29 | 96 | 104 | 104 |
| 5 | 82+ | 77+ | 104+ | 30 | 70 | 104 | 77 |
| 6 | 89 | 91 | 90 | 31 | 91+ | 70+ | 92+ |
| 7 | 39 | 45 | 50 | 32 | 103 | 69 | 91 |
| 8 | 93+ | 104+ | 103+ | 33 | 85+ | 72+ | 104+ |
| 9  | 104+  | 63+   | 104+  | 34 | 104+  | 104+  | 74+   |
| 10 | 81+   | 104+  | 69+   | 35 | 67    | 104   | 68    |
| 11 | 104+  | 104+  | 104+  | 36 | 104+  | 104+  | 104+  |
| 12 | 104+  | 83+   | 40+   | 37 | 87+   | 104+  | 104+  |
| 13 | 104+  | 104+  | 104+  | 38 | 89+   | 104+  | 104+  |
| 14 | 78+   | 104+  | 104+  | 39 | 104+  | 81+   | 64+   |
| 15 | 86    | 55    | 94    | 40 | 34    | 104   | 54    |
| 16 | 76+   | 87+   | 74+   | 41 | 103   | 73    | 84    |
| 17 | 102   | 104   | 80    | 42 | 80    | 104   | 73    |
| 18 | 45    | 79    | 104   | 43 | 94    | 104   | 104   |
| 19 | 104+  | 104+  | 104+  | 44 | 104+  | 101+  | 94+   |
| 20 | 76+   | 84+   | 78+   | 45 | 80    | 81    | 76    |
| 21 | 72    | 95    | 104   | 46 | 73    | 104   | 66    |
| 22 | 92    | 104   | 102   | 47 | 104+  | 98+   | 73+   |
| 23 | 55+   | 104+  | 104+  | 48 | 49+   | 83+   | 77+   |
| 24 | 89    | 104   | 104   | 49 | 88+   | 79+   | 99+   |
| 25 | 103   | 91    | 104   | 50 | 104+  | 104+  | 79+   |

+ denote right censored times

Table 2. Number of litters of various censoring status combinations

| Censoring Status (T, C1, C2) | Number of Litters |
|---|---|
| (0, 0, 0) | 1 |
| (0, 0, 1) | 3 |
| (0, 1, 0) | 6 |
| (1, 0, 0) | 2 |
| (0, 1, 1) | 11 |
| (1, 0, 1) | 2 |
| (1, 1, 0) | 2 |
| (1, 1, 1) | 23 |

Tumor occurrence is denoted as 0 and censored is denoted as 1

Table 3. Maximum likelihood estimates and standard errors.

| Parameter | Estimate | Standard Error |
|---|---|---|
| $\alpha$ | 0.900 | 0.101 |
| Scale (T) | 112.035 | 6.484 |
| Shape (T) | 4.393 | 0.879 |
| Scale (C1) | 157.641 | 28.649 |
| Shape (C1) | 3.568 | 1.150 |
| Scale (C2) | 154.119 | 25.667 |
| Shape (C2) | 2.890 | 0.792 |